\documentclass[11pt]{article}
\PassOptionsToPackage{obeyspaces}{url}
\usepackage[hidelinks]{hyperref}
\usepackage{bookmark}
\bookmarksetup{open,numbered,depth=5}
\usepackage{amsfonts}
\usepackage{amsmath}
\usepackage{amssymb, amsthm, enumitem}
\usepackage{mathrsfs}
\usepackage{breakurl}
\usepackage{pgf,tikz}
\usepackage{bm}

\usetikzlibrary{calc}

\usepackage{etoolbox} 
\patchcmd{\thebibliography}{\leftmargin\labelwidth}{\leftmargin\labelwidth\addtolength\itemsep{-0.1\baselineskip}}{}{}

\usepackage{mathtools}
\usepackage{xstring}

\author{Zijian Xu\thanks{School of Mathematical Sciences, Peking University, Beijing 100871. \texttt{2200010770@stu.pku.edu.cn}.}}

\title{On $k$-modal subsequences}
\date{}

\usepackage[nameinlink]{cleveref}

\newtheorem{theorem}{Theorem}
\newtheorem{lemma}[theorem]{Lemma}

\newtheorem{observation}[theorem]{Observation}
\newtheorem{proposition}[theorem]{Proposition}

\newcommand*{\eqdef}{\stackrel{\mbox{\normalfont\tiny def}}{=}} 

\newcommand*{\N}{\mathbb{N}}                                    
\newcommand*{\Z}{\mathbb{Z}}                                    
\newcommand*{\R}{\mathbb{R}}                                    

\newcommand*{\cS}{\mathcal{S}}

\newcommand*{\sN}{\mathsf{N}}

\newcommand*{\tS}{\widetilde{S}}
\newcommand*{\vH}{\mathscr{H}}

\newcommand*{\cO}{\mathcal{O}}
\newcommand*{\cC}{\mathcal{C}}
\newcommand*{\cI}{\mathcal{I}}
\newcommand*{\cD}{\mathcal{D}}
\newcommand*{\cF}{\mathcal{F}}

\newcommand*{\cL}{\mathcal{L}}

\newcommand*{\Gr}{\Gamma^Q_{\operatorname{rich}}}
\newcommand*{\Gp}{\Gamma^Q_{\operatorname{poor}}}

\newcommand*{\In}[1]{\{0,1,\dots,#1\}}

\crefname{enumi}{step}{steps}
\crefname{part}{part}{parts}

\allowdisplaybreaks[4]

\begin{document}

\maketitle

\begin{abstract}
    A $k$-modal sequence is a sequence of real numbers that can be partitioned into $k+1$ (possibly empty) monotone sections such that adjacent sections have opposite monotonicities. For every positive integer $k$, we prove that any sequence of $n$ pairwise distinct real numbers contains a $k$-modal subsequence of length at least $\sqrt{(2k+1)(n-\frac14)} - \frac{k}{2}$, which is tight in a strong sense. This confirms an old conjecture of F.R.K.~Chung (\emph{J.~Combin.~Theory Ser.~A}, 29(3):267–279, 1980). 
\end{abstract}

\section{Introduction} \label{sec:intro}

Given a sequence of real numbers, does there exist a long subsequence under some fixed structure $\cS$? This kind of Ramsey-style problems have attracted lots of attention over the years. Probably the first result of such a problem (in which $\cS = \text{``monotonicity’’}$) is the well-known Erd\H{o}s--Szekeres theorem \cite{erdos_szekeres}. It asserts that every sequence of length $n$ contains a monotone subsequence of length at least $\bigl\lceil\sqrt{n}\bigr\rceil$, a tight lower bound. Moreover, this monotone subsequence of a sequence problem can be naturally extended to analogy problem concerning the convex subset of a point set. Indeed, the convex subset problem was initiated in \cite{erdos_szekeres}, and was reiterated in \cite{chung_graham,kleitman_pachter,toth_valtr_1998,toth_valtr_2005,mojarrad_vlachos,norin_yuditsky} since then. This direction of study culminated in the result that any $2^{n+o(n)}$ points in $\R^2$ contain $n$ points in convex position, a breakthrough by Suk \cite{suk}. Other studies are around the case $\cS = \text{``higher order monotonicity’’}$ (\cite{elias_matousek,moshkovitz_shapira}) and the case $\cS = \text{``convexity together with an empty upper region’’}$ (\cite{valtr,cerny}). There are also recent results on monotone subarrays (\cite{bucic_sudakov_tran,ciceksiz_jin_raty_tomon}). 

Through out this paper, $k$ stands for some nonnegative integer. We work on the aforementioned problem in which $\cS = \text{``piecewise monotonicity’’}$. To be specific, a \emph{generic sequence of length $n$} is a sequence of distinct $a_1, \dots, a_n \in \R$, and a \emph{$k$-modal subsequence of length $t$} of a generic sequence is a subsequence $b_1, \dots, b_t$ of $a_1, \dots, a_n$ such that there exist $1 = i_0 \le i_1 \le \dots \le i_k \le i_{k+1} = t$ with
\begin{itemize}
    \item $b_{i_j} < b_{i_j+1} < \dots < b_{i_{j+1}}$ if $j$ is odd and $b_{i_j} > b_{i_j+1} > \dots > b_{i_{j+1}}$ if $j$ is even, or
    \item $b_{i_j} > b_{i_j+1} > \dots > b_{i_{j+1}}$ if $j$ is odd and $b_{i_j} < b_{i_j+1} < \dots < b_{i_{j+1}}$ if $j$ is even. 
\end{itemize}
In particular, a unimodal subsequence refers to a $1$-modal subsequence. 

This definition can be viewed as a natural extension of the Erd\H{o}s--Szekeres monotone subsequence problem. By definition, a monotone subsequence is a $0$-modal subsequence. Let $\rho(n; k)$ be the largest integer such that any generic sequence of length $n$ contains a $k$-modal subsequence of length $\rho(n; k)$. The Erd\H{o}s--Szekeres theorem asserts that $\rho(n; 0) = \bigl\lceil \sqrt{n} \bigr\rceil$. Chung \cite{chung} resolved the unimodal case by proving $\rho(n; 1) = \Bigl\lceil \sqrt{3n-\frac{3}{4}} - \frac{1}{2}\Bigr\rceil$. A related probabilistic problem was studied by Steele \cite{steele}. 

Chung \cite{chung} also mentioned that $\rho(n; k) \le \sqrt{(2k+1)n}$, and conjectured that this is asymptotically tight when $k$ is fixed and $n$ goes to infinity. By studying $k$-modal subsequences which are increasing from the beginning to the first modal, Gong et al.~\cite{bukh_dong_gong_gu_yang_yu} established $\rho(n;k) \ge \sqrt{2kn}$. In this paper, we resolve Chung's longstanding conjecture by proving the following pair of theorems: 

\begin{theorem} \label{thm:modal_subseq}
    For any integers $k \ge 0$ and $n \ge 1$, we have $\rho(n; k) \ge \Bigl\lceil \sqrt{(2k+1)(n-\frac14)} - \frac{k}{2} \Bigr\rceil$. 
\end{theorem}

\begin{theorem} \label{thm:modal_construction}
    For any integers $k \ge 1$ and $n \ge 10k^3$, we have $\rho(n;k) \le \Bigl\lceil \sqrt{(2k+1)(n-\frac14)} - \frac{k}{2} \Bigr\rceil + 1$. 
\end{theorem}

The $n \ge 10k^3$ threshold in \Cref{thm:modal_construction} is not optimal, but one cannot expect the same bound for small $n$. Indeed, it is obvious from the definition that $\rho(n;k) = n$ when $n \le k$. Also, it might very well be possible to remove the error term ``$+1$'' in \Cref{thm:modal_construction} through a more careful analysis of our explicit construction, but we decide not to put into that level of effort. 

In the cases when $k = 0$ or $1$, \Cref{thm:modal_subseq} is best possible as it matches the aforementioned values $\rho(n;0)$ and $\rho(n;1)$. For larger $k$, \Cref{thm:modal_construction} shows that \Cref{thm:modal_subseq} is tight up to a small additive constant, which confirms Chung's conjectured $\rho(n;k) = \bigl( 1 + o(1) \bigr) \sqrt{(2k+1)n}$ in a strong sense. 

Probably the most famous proof of the lower bound $\rho(n;0) \ge \sqrt{n}$ is done by the labeling trick introduced by Seidenberg \cite{seidenberg}, where each term $a_i$ of $(a_i)_{i=1}^n$ is associated with the label
\[
\begin{pmatrix}
    x_i \\ y_i
\end{pmatrix}
\eqdef 
\begin{pmatrix}
    \text{the length of the longest increasing subsequence ending at $a_i$} \\
    \text{the length of the longest decreasing subsequence ending at $a_i$}
\end{pmatrix}. 
\]
Then the lower bound follows immediately by observing that the labels are pairwise distinct. The lower bounds on $\rho(n;1)$ due to Chung and $\rho(n;k)$ due to Gong et al.~are also established by more sophisticated labeling arguments. A main novelty of this work is that our proof of \Cref{thm:modal_subseq} does not utilize any kind of labeling. Inspired by the notable theorem due to Dilworth \cite{dilworth}, we shall prove \Cref{thm:modal_subseq} by analyzing the underlying poset structures behind $k$-modal subsequences. It might be of independent interest to figure out a proof of \Cref{thm:modal_subseq} via some kind of labels. 

\smallskip

Instead of proving \Cref{thm:modal_subseq,thm:modal_construction} directly, we consider the geometric generalization, which turns out to be easier to work with. For every point $p = (a, b)$ on the Euclidean plane $\R^2$, denote $x(p) \eqdef a$ and $y(p) \eqdef b$. For any pair of points $p, q \in \R^2$, we introduce two partial orders: 
\begin{itemize}
    \item The ``northeast'' one $\prec_{_{\text{NE}}}$ satisfies $p \prec_{_{\text{NE}}} q$ if and only if $x(p) < x(q)$ and $y(p) < y(q)$. 
    \item The ``southeast'' one $\prec_{_{\text{SE}}}$ satisfies $p \prec_{_{\text{SE}}} q$ if and only if $x(p) < x(q)$ and $y(p) > y(q)$. 
\end{itemize}
Among finite subsets of $\R^2$, an \emph{increasing path} is a chain under $\prec_{_{\text{NE}}}$, a \emph{decreasing path} is a chain under $\prec_{_{\text{SE}}}$, and a \emph{monotone path} is either an increasing path or a decreasing path. The number of points in a monotone path is referred to as its \emph{length}. For instance, the empty set $\varnothing$ is a monotone path of length $0$, and the set $\{(1, 5), (2, 4), (8, \pi)\}$ is a decreasing path of length $3$. 

Let $S \subset \R^2$ be a finite point set. We call $S$ \emph{generic} if $x(p) \neq x(q)$ and $y(p) \neq y(q)$ for any pair of distinct points $p, q \in S$. The $i$-th point of $S$ refers to the element of $S$ with the $i$-th smallest $x$-coordinate. A $(k+1)$-partition of $S$ is a $(k+1)$-tuple $(S_0, \dots, S_k)$ with $S = S_0 \cup \dots \cup S_k$ such that $\max\limits_{p \in S_i} x(p) < \min\limits_{q \in S_j} x(q)$ whenever $0 \le i < j \le k$. Here each $S_i \, (0 \le i \le k)$ is possibly empty, and we apply the conventions $\max\limits_{p \in \varnothing} x(p) = -\infty, \, \min\limits_{q \in \varnothing} x(q) = +\infty$. 

\smallskip
A generic $P \subset \R^2$ is a $k$-\emph{modal path} if $P$ admits a $(k+1)$-partition $(P_0, \dots, P_k)$ such that
\begin{itemize}
    \item for $i = 0, \dots, k$, each section $P_i$ is a monotone path, and
    \item for $i = 1, \dots, k$, the monotonicity of $P_{i-1}, P_i$ are opposite. 
\end{itemize}
For a $k$-modal path $P$, we distinguish two possible classes of it: 
\begin{itemize}
    \item If $P_0, P_2, \dots$ are increasing paths and $P_1, P_3, \dots$ are decreasing paths, then $P$ is a $(+k)$-path. 
    \item If $P_0, P_2, \dots$ are decreasing paths and $P_1, P_3, \dots$ are increasing paths, then $P$ is a $(-k)$-path. 
\end{itemize}
Note that the two classes are not disjoint as each section $P_i$ might have size at most $1$. The \emph{length} of a $k$-modal path refers to its cardinality as a finite set. 

By considering the generic point set $\{(i, a_i) : i = 1, \dots, n\} \subset \R^2$ generated by a generic sequence $(a_i)_{i=1}^n$, and the generic sequence $(y_i)_{i=1}^n$ induced by a generic point set $\{(x_i, y_i) : i = 1, \dots, n\} \subset \R^2$, it is easily seen that \Cref{thm:modal_subseq,thm:modal_construction} are equivalent to 
\Cref{thm:modal_subpth,thm:modal_setconstruct} below, respectively. 

\begin{theorem} \label{thm:modal_subpth}
    If nonempty $S \subset \R^2$ is generic, then $S$ contains a $k$-modal path subset $P$ of length
    \[
    |P| \geq \Bigl\lceil \sqrt{(2k+1)(|S|-\tfrac14)} - \tfrac{k}{2} \Bigr\rceil. 
    \]
\end{theorem}

\begin{theorem} \label{thm:modal_setconstruct}
    Given integers $k \ge 1$ and $n \ge 10k^3$. There exists an $n$-point generic $S \subset \R^2$ such that 
    \[
    |P| \leq \Bigl\lceil \sqrt{(2k+1)(|S|-\tfrac14)} - \tfrac{k}{2} \Bigr\rceil + 1 
    \]
    holds for the length of any $k$-modal path subset $P$ of $S$. 
\end{theorem}

The rest of this paper is organized as follows. \Cref{sec:finecover,sec:prop} (also \Cref{sec:proof_inj}) are devoted to the proof of \Cref{thm:modal_subpth}. In \Cref{sec:finecover}, we introduce \emph{fine coverings}, the main tool in our proof of \Cref{thm:modal_subpth}. We also prove \Cref{thm:modal_subpth} in \Cref{sec:finecover} assuming a key property (\Cref{thm:fine_covering}) of fine coverings, and then establish \Cref{thm:fine_covering} in \Cref{sec:prop} (the proof of a technical lemma is postponed to \Cref{sec:proof_inj}). Finally, we present our constructive proof of \Cref{thm:modal_setconstruct} in \Cref{sec:construction}. 

\section{Modal paths and fine coverings} \label{sec:finecover}

Observe that antichains with respect to ``$\prec_{_\text{NE}}$'' are exactly chains with respect to ``$\prec_{_\text{SE}}$'', and vice versa. We thus deduce from Dilworth's theorem (see \cite{dilworth}) the following result. 

\begin{lemma} \label{lem:dilworth_of_paths}
    Let $S \subset \R^2$ be a generic point set. Then the maximum length of increasing paths of $S$ is equal to the minimum number of decreasing paths needed to cover $S$. 
\end{lemma}

Let $S \subset \R^2$ be generic. A \emph{fine covering} of $S$ by $(+k)$-paths is a multiset $\cC$ of $(+k)$-paths such that the $i$-th sections of all $k$-modal paths in $\cC$ form a partition of $S$ for each $i \in \In{k}$. The \emph{size} of $\cC$ is the cardinality of $\cC$ as a multiset, denoted by $|\cC|$. Such a covering will be abbreviated as a $(+k)$-\emph{covering}. The reason we think of $\cC$ as a multiset is that two $(+k)$-paths in $\cC$ may be the same as a set but admit different partitions into sections. For instance, if $S = P = Q = \{(1, 1), (2, 2)\}$ with $P_0 = \varnothing, \, P_1 = S$ and $Q_0 = S, \, Q_1 = \varnothing$, then $\cC = \{P, Q\}$ is a $(+1)$-covering of $S$. We also define $(-k)$-\emph{coverings}, fine coverings by $(-k)$-paths, in a similar way. 

For any generic point set $S \subset \R^2$, we introduce the following definitions: 
\begin{itemize}
    \item Denote by $M_+(S), M_-(S)$ the maximum length of $(+k)$- and $(-k)$-paths of $S$, respectively. 
    \item Denote by $m_+(S), m_-(S)$ the minimum size of $(+k)$- and $(-k)$-coverings of $S$, respectively. 
\end{itemize}

\begin{proposition} \label{prop:part_1}
     Let $S$ be a generic point set with $|S| = n$. Then we have
     \begin{itemize}
         \item $m_+(S) \leq M_+(S) + k$ and $m_-(S) \leq M_-(S) + k$ if $k$ is odd, and
         \item $m_+(S) \leq M_-(S) + k$ and $m_-(S) \leq M_+(S) + k$ if $k$ is even. 
     \end{itemize}
\end{proposition}

\begin{proof}
We are going to construct a generic point set and apply \Cref{lem:dilworth_of_paths}. Without loss of generality, we assume that $S\subseteq \R\times (0,1)$ and $|x(p) - x(q)| \geq k+1$ for any pair of distinct points $p, q \in S$. 

We inductively define generic point sets $S_0, \dots, S_k$ as follows: Set $S_0 \eqdef S$. For each $i = 1, \dots, k$, assume $S_{i-1}$ is defined and set $S_i \eqdef r_i(S_{i-1})$. Here $r_i$, a composition of the reflection over the line $y = i$ and the shifting under the vector $(1, 0)$, sends $(x, y)$ to $(x+1, 2i-y)$. The trajectory of each point of $S$ induces a natural bijection $\tau_i \colon S_0 \to S_i$. Indeed, the $i$-th point of $S_0$ is mapped to the $i$-th point of $S_i$ under $\tau_i$. \hyperlink{figone}{Figure~1} illustrates this process, in which $|S| = 5$ and $k = 3$. 

\begin{center}
\begin{tikzpicture}[x = 0.3cm, y = 0.3cm]
    \clip (-1, -1) rectangle (24, 3);
    \draw[ultra thin] (0, 0) rectangle (19, 1);
    \filldraw[color = black] (0.5, 0.4) circle (0.075);
    \filldraw[color = black] (4.8, 0.7) circle (0.075);
    \filldraw[color = black] (9.0, 0.2) circle (0.075);
    \filldraw[color = black] (14.1, 0.8) circle (0.075);
    \filldraw[color = black] (18.5, 0.3) circle (0.075);
    \node at (19.8, 0.5) {\begin{footnotesize}$S_0$\end{footnotesize}};
\end{tikzpicture}
\begin{tikzpicture}[x = 0.3cm, y = 0.3cm]
    \clip (-1, -1) rectangle (24, 3);
    \draw[ultra thin] (0, 0) rectangle (19, 1);
    \filldraw[color = black] (0.5, 0.4) circle (0.075);
    \filldraw[color = black] (4.8, 0.7) circle (0.075);
    \filldraw[color = black] (9.0, 0.2) circle (0.075);
    \filldraw[color = black] (14.1, 0.8) circle (0.075);
    \filldraw[color = black] (18.5, 0.3) circle (0.075);
    \node at (19.8, 0.5) {\begin{footnotesize}$S_0$\end{footnotesize}};
    \draw[ultra thin] (1, 1) rectangle (20, 2);
    \filldraw[color = black] (1.5, 1.6) circle (0.075);
    \filldraw[color = black] (5.8, 1.3) circle (0.075);
    \filldraw[color = black] (10.0, 1.8) circle (0.075);
    \filldraw[color = black] (15.1, 1.2) circle (0.075);
    \filldraw[color = black] (19.5, 1.7) circle (0.075);
    \node at (20.8, 1.5) {\begin{footnotesize}$S_1$\end{footnotesize}};
\end{tikzpicture}
\begin{tikzpicture}[x = 0.3cm, y = 0.3cm]
    \clip (-1, -1) rectangle (24, 5);
    \draw[ultra thin] (0, 0) rectangle (19, 1);
    \filldraw[color = black] (0.5, 0.4) circle (0.075);
    \filldraw[color = black] (4.8, 0.7) circle (0.075);
    \filldraw[color = black] (9.0, 0.2) circle (0.075);
    \filldraw[color = black] (14.1, 0.8) circle (0.075);
    \filldraw[color = black] (18.5, 0.3) circle (0.075);
    \node at (19.8, 0.5) {\begin{footnotesize}$S_0$\end{footnotesize}};
    \draw[ultra thin] (1, 1) rectangle (20, 2);
    \filldraw[color = black] (1.5, 1.6) circle (0.075);
    \filldraw[color = black] (5.8, 1.3) circle (0.075);
    \filldraw[color = black] (10.0, 1.8) circle (0.075);
    \filldraw[color = black] (15.1, 1.2) circle (0.075);
    \filldraw[color = black] (19.5, 1.7) circle (0.075);
    \node at (20.8, 1.5) {\begin{footnotesize}$S_1$\end{footnotesize}};
    \draw[ultra thin] (2, 2) rectangle (21, 3);
    \filldraw[color = black] (2.5, 2.4) circle (0.075);
    \filldraw[color = black] (6.8, 2.7) circle (0.075);
    \filldraw[color = black] (11.0, 2.2) circle (0.075);
    \filldraw[color = black] (16.1, 2.8) circle (0.075);
    \filldraw[color = black] (20.5, 2.3) circle (0.075);
    \node at (21.8, 2.5) {\begin{footnotesize}$S_2$\end{footnotesize}};
\end{tikzpicture}
\begin{tikzpicture}[x = 0.3cm, y = 0.3cm]
    \clip (-1, -1) rectangle (24, 5);
    \draw[ultra thin] (0, 0) rectangle (19, 1);
    \filldraw[color = black] (0.5, 0.4) circle (0.075);
    \filldraw[color = black] (4.8, 0.7) circle (0.075);
    \filldraw[color = black] (9.0, 0.2) circle (0.075);
    \filldraw[color = black] (14.1, 0.8) circle (0.075);
    \filldraw[color = black] (18.5, 0.3) circle (0.075);
    \node at (19.8, 0.5) {\begin{footnotesize}$S_0$\end{footnotesize}};
    \draw[ultra thin] (1, 1) rectangle (20, 2);
    \filldraw[color = black] (1.5, 1.6) circle (0.075);
    \filldraw[color = black] (5.8, 1.3) circle (0.075);
    \filldraw[color = black] (10.0, 1.8) circle (0.075);
    \filldraw[color = black] (15.1, 1.2) circle (0.075);
    \filldraw[color = black] (19.5, 1.7) circle (0.075);
    \node at (20.8, 1.5) {\begin{footnotesize}$S_1$\end{footnotesize}};
    \draw[ultra thin] (2, 2) rectangle (21, 3);
    \filldraw[color = black] (2.5, 2.4) circle (0.075);
    \filldraw[color = black] (6.8, 2.7) circle (0.075);
    \filldraw[color = black] (11.0, 2.2) circle (0.075);
    \filldraw[color = black] (16.1, 2.8) circle (0.075);
    \filldraw[color = black] (20.5, 2.3) circle (0.075);
    \node at (21.8, 2.5) {\begin{footnotesize}$S_2$\end{footnotesize}};
    \draw[ultra thin] (3, 3) rectangle (22, 4);
    \filldraw[color = black] (3.5, 3.6) circle (0.075);
    \filldraw[color = black] (7.8, 3.3) circle (0.075);
    \filldraw[color = black] (12.0, 3.8) circle (0.075);
    \filldraw[color = black] (17.1, 3.2) circle (0.075);
    \filldraw[color = black] (21.5, 3.7) circle (0.075);
    \node at (22.8, 3.5) {\begin{footnotesize}$S_3$\end{footnotesize}};
\end{tikzpicture}
\begin{tikzpicture}
    \clip (-8, -0.25) rectangle (8, 0.5);
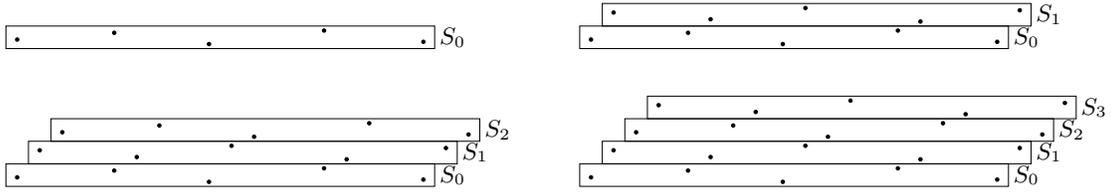
    \node at (0, 0) {\textbf{\hypertarget{figone}{Figure 1:}} An evolution from $S_0$ to $S_3$. };
\end{tikzpicture}
\end{center}

Let $\tS \eqdef \bigcup\limits_{i=0}^k S_i$. Obviously, $\tS$ is generic. We introduce a map $g$ such that
\[
\{\text{decreasing paths of $\tS$}\} \,\, \overset{g}{\longrightarrow} \,\, \{\text{$k$-modal paths of $S$}\}. 
\]
To be specific, for any decreasing path $P \subseteq \widetilde{S}$, we define $g(P) \eqdef \bigcup\limits_{i=0}^k \bigl(\tau_i^{-1} (P \cap S_i) \bigr) \subseteq S$ (note that $\tau_0$ is the identity map on $S_0$ and $\tau(U) = \{\tau(u) : u \in U\}$). Moreover, as a $k$-modal path, the $i$-th section of $g(P)$ is $\tau_{k-i}^{-1}(P \cap S_{k-i})$ for $i = 0, 1, \dots, k$. Observe that $g(P)$ is a $(+k)$-path when $k$ is odd, and $g(P)$ is a $(-k)$-path when $k$ is even. (Also, $|g(P)| = |P|$ because the trajectory $\{p, \tau_1(p), \dots, \tau_k(p)\}$ is an increasing path for any $p \in P$, but we shall not use this property.)  

For any increasing path $P$ of $\widetilde{S}$, similar to the definition of $g$, we associate with it a $(+k)$-path $T_P \eqdef \bigcup\limits_{i=0}^k \bigl( \tau_i^{-1}(P \cap S_i) \bigr)$, whose $i$-th section is given by $\tau_i^{-1}(P \cap S_i)$. Notice that $T_P$ is always a $(+k)$- rather than a $(-k)$-path because its $0$-th section $\tau_0^{-1}(P \cap S_0) = P \cap S_0$ is increasing. 

We claim that $|T_P| \ge |P| - k$. To see this, it suffices to show for $i = 1, \dots, k$ that
\[
\Biggl|\tau_i^{-1}(P \cap S_i) \cap \biggl( \bigcup_{j=0}^{i-1} \bigl( \tau_j^{-1}(P \cap S_j) \bigr) \biggr) \Biggr| \leq 1. 
\]
Indeed, this is true because the intersection happens only possibly at the last point of $\tau_i^{-1}(P \cap S_i)$ and the first point of $\bigcup\limits_{j=0}^{i-1} \bigl( \tau_j^{-1}(P \cap S_j) \bigr)$. Thus, the claim holds. 

Since the length of each $(+k)$-path of $S$ is upper bounded by $M_+(S)$, the claimed $|T_P| \ge |P| - k$ implies that any increasing path on $\widetilde{S}$ is of length at most $M_+(S) + k$. We thus obtain from \Cref{lem:dilworth_of_paths} the existence of $\cF$, a covering of $\tS$ by $M_+(S) + k$ many (possibly empty) decreasing paths. Assume further that each point of $S$ is covered exactly once in $\cF$. Then $g$ induces a covering $\cF'\eqdef g(\cF)$ of $S$ by $M_+(S) + k$ many $(\pm k)$-paths (``$\pm$'' depends on the parity of $k$). So, $M_+(S)+k \geq m_\pm(S)$. 

Suppose $k$ is odd. Then $\cF'$ gives a $(+k)$-covering of $S$ with $|\cF'| \le M_+(S) + k$. This implies that $m_+(S) \le M_+(S) + k$. Let $\sigma$ be the reflection over the line $y = \frac{1}{2}$ in the plane. By applying the same arguments to $\sigma(S)$ we deduce that $m_-(S) \le M_-(S) + k$. The case when $k$ is even can be proved similarly, and the proof of \Cref{prop:part_1} is complete. 
\end{proof}

For any generic point set $S \subset \R^2$, write $m(S) \eqdef \max \bigl\{ m_+(S), m_-(S) \bigr\}$. In other words, there are
\begin{itemize}
    \item a $(+k)$-covering of $S$ by at most $m(S)$ many $(+k)$-paths, and
    \item a $(-k)$-covering of $S$ by at most $m(S)$ many $(-k)$-paths. 
\end{itemize}

\begin{theorem} \label{thm:fine_covering}
    For any generic $S$, we have $m(S) \geq \sqrt{(2k+1)(|S|-\frac14)} + \frac{k}{2}$. 
\end{theorem}

\begin{proof}[Proof of \Cref{thm:modal_subpth} assuming \Cref{thm:fine_covering}]
    Let $P$ be a $k$-modal path of $S$ whose length is maximized. Then \Cref{prop:part_1} and \Cref{thm:fine_covering} implies that $|P| \ge m(S) - k \ge \sqrt{(2k+1)(|S|-\frac14)} - \frac{k}{2}$. Since $|P|$ takes integer values, we conclude that $|P| \ge \Bigl\lceil \sqrt{(2k+1)(|S|-\frac14)} - \frac{k}{2} \Bigr\rceil$. 
\end{proof}


\section{Proof of \texorpdfstring{\Cref{thm:fine_covering}}{Theorem~\ref{thm:fine_covering}}} \label{sec:prop}

The definition of $m(S)$ implies that there exists a $(+k)$-covering $\cC_+$ and a $(-k)$-covering $\cC_-$ of $S$ with $|\cC_+| \le m(S)$ and $|\cC_-| \le m(S)$. One extra assumption on $\cC_+, \cC_-$ will be imposed later. 

Recall that every $k$-modal path admits some $(k+1)$-partition in which each section is of some fixed monotonicity. Throughout this section, whenever we mention a $k$-modal path $\bullet$, we implicitly assume that $\bullet$ has the $(k+1)$-partition $\bullet_0 \cup \dots \cup \bullet_k$. The reason we write $\bullet$ rather than $P$ here is that $\bullet$ could be $P, Q, R, L$ with possible superscripts later in the text. 

Assume without loss of generality that $\cC_+ \cup \cC_-$ does not contain the empty $k$-modal path. Due to the existence of a $k$-modal path that is both a $(+k)$- and a $(-k)$-path (such as a $k$-point generic set), we denote by $\cC_+ \cup \cC_-$ the union of multisets (which implies that $|\cC_+ \cup \cC_-| = |\cC_+| + |\cC_-|$). 

For each $P \in \cC_+ \cup \cC_-$, define a map $\sN_P \colon P \setminus P_k \to \cC_+ \cup \cC_-$ such that: 
\begin{itemize}
    \item If $p \in P_i$ and $P \in \cC_+$, then $\sN_P(p)$ is the unique $Q \in \cC_-$ satisfying $p \in Q_{i+1}$. 
    \item If $p \in P_i$ and $P \in \cC_-$, then $\sN_P(p)$ is the unique $Q \in \cC_+$ satisfying $p \in Q_{i+1}$. 
\end{itemize}
We remark that the uniqueness follows from the definition of fine coverings. 

Write $s(P)$ as the first point of $P$, and denote by $i(P)$ the unique index $i$ such that $s(P) \in P_i$. Define $\cD \eqdef \{P \in \cC_+ \cup \cC_- : i(P) < k\}$. Let $\varphi \colon \cD \to \cC_+ \cup \cC_-$ be the map given by $\varphi(P) \eqdef \sN_P\bigl(s(P)\bigr)$. 

We claim that the action of $\varphi$ is acyclic. Formally, there exists no $\ell \in \N_+$ and $k$-modal path $P$ such that the $\ell$-th iteration $\varphi^{(\ell)}(P) = P$. To see this, write $Q \eqdef \varphi(P)$. Then $s(P) \in Q$, and hence $x\bigl(s(Q)\bigr) \le x\bigl(s(P)\bigr)$ with $x\bigl(s(Q)\bigr) = x\bigl(s(P)\bigr)$ implying $s(Q) = s(P)$ and so $i(Q) = i(P) + 1$. Thus, 
\[
\Bigl( x\bigl(s(Q)\bigr), -i(Q) \Bigr) \prec \Bigl( x\bigl(s(P)\bigr), -i(P) \Bigr), 
\]
where ``$\prec$'' stands for the canonical lexicographical order on $\Z^2$. This verifies the claim. 

Let $\varphi^{(t)}$ be the $t$-th iteration of $\varphi$, where $\varphi^{(1)} = \varphi$ and $\varphi^{(0)}$ is the identity. For each $P \in \cC_+ \cup \cC_-$, we define its \emph{depth} $d(P)$ as follows: If $P \in (\cC_+ \cup \cC_-) \setminus \cD$, then $d(P) \eqdef 0$. Otherwise, when $P \in \cD$, set $d(P) \eqdef \ell$ be the unique positive integer such that $\varphi^{(\ell-1)}(P) \in \cD$ and $\varphi^{(\ell)}(P) \in (\cC_+ \cup \cC_-) \setminus \cD$. It is worth mentioning that $\ell$ exists because the action of $\varphi$ is acyclic (the previous claim). 

\begin{observation} \label{obs:depth_0th_section}
    For any $P \in \cC_+ \cup \cC_-$, if $d(P) < k$, then $P_0 = \varnothing$. 
\end{observation}

\begin{proof}
    It follows from the definitions of $\varphi$ and $d$ that $i(P) = k$ when $d(P) = 0$. Then it is easy to show (via induction on $j$) that $i(P) \geq k-j$ when $d(P) = j$, for any nonnegative integer $j$. 
\end{proof}

We remark that the converse of \Cref{obs:depth_0th_section} does not hold, as $d(P) \ge k$ does not necessarily imply that $P_0 \neq \varnothing$. Notice that the $\sN_P, \cD, \varphi, d$ defined in this section all depend on $\cC_+, \cC_-$. 

\begin{lemma} \label{lem:inj}
    We can choose $\cC_+, \cC_-$ appropriately so that $\varphi$ is injective. 
\end{lemma}

The proof is subtle and is included in \Cref{sec:proof_inj}. Thanks to \Cref{lem:inj}, we may assume that $\varphi$ is injective. This is the extra assumption on $\cC_+, \cC_-$ that we claimed at the beginning of the section. 

Since $\varphi$ is injective, the action of $\varphi$ results in a finite number of disjoint chain-like orbits. Let $\cO_1, \dots, \cO_r$ be all such orbits containing at least $k+1$ elements ($k$-modal paths). For each of these orbits $\cO$, assume $\cO = \{P^0, P^1, \dots, P^t\}$ with $t \ge k$ and $d(P^i) = i$ for $i = 0, 1, \dots, t$. That is, 
\[
P^0 = \varphi(P^1) = \dots = \varphi^{(t)}(P^t) \quad \text{and} \quad \varphi(P^0) \in (\cC_+ \cup \cC_-) \setminus \cD. 
\]
Associate with $\cO$ the modal path $P^{\cO}$ whose $(i-k)$-th section is given by $P^{i-k}_0 \, (i = k, k+1, \dots, t)$. Indeed, $P^{\cO}$ is a $(t-k)$-modal path, and we define $P^{\cO} \eqdef \varnothing$ when $k > t$. To see that the associated $P^{\cO}$ is well-defined, we must verify for all $k \le i \le t$ that the monotonicity of $P^i_0, P^{i+1}_0$ are opposite, and for all $k \le i < j \le t$ that $\max x(P_0^i) < \min x(P_0^j)$. Notice that we need to consider not only the consecutive pairs because each $P^{\ell}_0$ could possibly be empty. From $\varphi(P^{\ell+1}) = P^{\ell}$ we deduce that
\begin{itemize}
    \item $\min x(P^j_0) \geq x\bigl(s(P^j)\bigr) \geq x\bigl(s(P^{j-1})\bigr) \geq \dots \geq x\bigl(s(P^{i+1})\bigr) \geq \max x(P^i_0)$, and
    \item $(P^i, P^{i+1}) \in \cC_+ \times \cC_-$ or $\cC_- \times \cC_+$, so $P^i_0, P^{i+1}_0$ are of opposite monotonicity. 
\end{itemize}
We shall write $R^i \eqdef P^{\cO_i} \, (i = 1, \dots, r)$ for brevity. Let $t_i+1 \ge k+1$ be the number of $k$-modal paths in the orbit $\cO_i$. Set $\ell_i \eqdef t_i - k + 1$, and so $R^i$ is an $(\ell_i - 1)$-modal path. 

\begin{observation} \label{obs:R_0th_section}
    The sections of $R^i$ consists of all $0$-th sections of $P \in \cC_+ \cup \cC_-$ with $d(P) \ge k$. 
\end{observation}

\begin{proof}
    This is straightforward from our definition of $P^{\cO}$. 
\end{proof}

\begin{observation} \label{obs:orbit_sum}
    $kr + (\ell_1+\dots+\ell_r) \le 2m(S)$. 
\end{observation}

\begin{proof}
From the definition of $\ell_i$ we obtain $k + \ell_i = t_i + 1 = |\cO_i|$. As disjoint orbits, $\cO_1, \dots, \cO_r$ form a partition of some submultiset of the multiset $\cC_+ \cup \cC_-$. It follows that
\[
kr + (\ell_1+\dots+\ell_r) = \sum_{i=1}^r (k+\ell_i) = \sum_{i=1}^r |\cO_i| \le |\cC_+ \cup \cC_-| = |\cC_+| + |\cC_-| \le 2m(S). \qedhere
\]
\end{proof}

For $i = 1, \dots, r$, we say that $R^i$ is a $+$-\emph{path} if $R^i$ is a $\bigl(+(\ell_i-1)\bigr)$-path, and $R^i$ is a $-$-\emph{path} if it is a $\bigl(-(\ell_i-1)\bigr)$-path. This $\pm$-path classification will be referred to as the \emph{sign} of $R^i$. Call
\begin{itemize}
    \item the sections of even indices of $+$-paths and of odd indices of $-$-paths \emph{increasing sections}, and
    \item the sections of odd indices of $+$-paths and of even indices of $-$-paths \emph{decreasing sections}. 
\end{itemize}
Let $\cI$ and $\cD$ be the multisets of all increasing and decreasing sections of $R^1, \dots, R^r$, respectively. 

Recall that a $(+0)$-covering of $S$ is a fine covering of $S$ by increasing paths, and a $(-0)$-covering of $S$ is a fine covering of $S$ by decreasing paths. 

\begin{proposition} \label{prop:covering_property}
    The set $\cI$ is a $(+0)$-covering of $S$, and the set $\cD$ is a $(-0)$-covering of $S$. 
\end{proposition}

\begin{proof}
    Since $\cC_+$ is a $(+k)$-covering of $S$, the collection of all $0$-th sections of $k$-modal paths in $\cC_+$ forms a $(+0)$-covering of $S$. It then suffices to verify that $\cI$ consists of all nonempty $0$-th sections of $P \in \cC_+$ and a bunch of empty sets. Indeed, from \Cref{obs:R_0th_section} we deduce that $\cI$ consists of all nonempty $0$-th sections of $P \in \cC_+$ with $d(P) \ge k$, and \Cref{obs:depth_0th_section} tells us that every $P \in \cC_+$ with $d(P) < k$ has an empty $0$-th section. We thus conclude that $\cI$ forms a $(+0)$-covering of $S$. Similarly, the set $\cD$ forms a $(-0)$-covering of $S$. 
\end{proof}

Let $K \eqdef \{(i, j) \in \N_+ : 1 \le i < j \le r\}$. We associate with each $p \in S$ a pair $\chi(p) \eqdef (i, j) \in K$ such that $p$ is the intersection of an increasing section and a decreasing section of $R^i$ and $R^j$. The existence of $R^i, R^j$ follows from \Cref{prop:covering_property}. Since we assumed $i < j$, we do not specify which section is from $R^i$ and which is from $R^j$. For $(i, j) \in K$, write $S_{i, j} \eqdef \{p \in S : \chi(p) = (i, j)\}$. Obviously, as $(i, j)$ goes through $K$, the sets $S_{i, j}$ form a partition of $S$. Define
\[
\ell(i, j) \eqdef \begin{cases}
    \bigl\lfloor \frac{\ell_i+\ell_j}{2} \bigr\rfloor &\text{if $R^i$ and $R^j$ are of opposite signs,} \\
    \bigl\lfloor \frac{\ell_i+\ell_j-1}{2} \bigr\rfloor &\text{if $R^i$ and $R^j$ are of the same sign.} 
\end{cases}
\]

\begin{proposition} \label{prop:Sij_bound}
    For each $(i, j) \in K$, we have $|S_{i, j}| \le \ell(i, j)$. 
\end{proposition}

\begin{proof}
    Let $p \in S_{i, j}$ and assume that $u, v$ are the unique pair of indices such that $p = R^i_u \cap R^j_v$. We claim that the parity of $u + v$ depends only on $(i, j)$. Indeed, from the definition we deduce that
    \[
    (u + v) \bmod 2 = \begin{cases}
        0 &\text{if $R^i$ and $R^j$ are of opposite signs,} \\
        1 &\text{if $R^i$ and $R^j$ are of the same sign.} 
    \end{cases}
    \]
    
    Let $p' \in S_{i, j}$ be a point other than $p$ and assume $p' = R^i_{u'} \cap R^j_{v'}$. We claim that $u + v \ne u' + v'$. To see this, we argue indirectly. Suppose $u+v = u'+v'$ and $(u', v') = (u+w, v-w)$. Due to their opposite monotonicity, $R^i_u$ and $R^j_v$ intersect in at most one point, and so $w \neq 0$. Assume without loss of generality that $w > 0$. This implies the following estimate on $x$-coordinates: 
    \[
    x(p') \le \max x(R^j_{v'}) < \min x(R^j_v) \le x(p) \le \max x(R^i_u) < \min x(R^i_{u'}) \le x(p'), 
    \]
    which is an obvious contradiction. The claim is thus proved. 

    Since $R^i$ is an $(\ell_i-1)$-modal path and $R^j$ is an $(\ell_j-1)$-modal path, by combining the analysis above we see that $|S_{i, j}|$ is upper bounded by either the total number of odds or the total number of evens among $\{0, 1, \dots, \ell_i + \ell_j - 2\}$, depending on the signs of $R^i, R^j$. The proof is done. 
\end{proof}

It is worth mentioning that our analysis in \Cref{prop:Sij_bound} is quite delicate. Indeed, we cannot prove it simply by establishing $|S_{i, j}| \le \min\{\ell_i, \ell_j\} \le \ell(i, j)$ as $|S_{i, j}| \le \min\{\ell_i, \ell_j\}$ does not hold in general. Also, the crude bound $|S_{i, j}| \le \max\{\ell_i, \ell_j\}$ is not enough for our proof of \Cref{thm:fine_covering}. 

\begin{proposition} \label{prop:lij_bound}
    We have
    \[
    \sum_{(i, j) \in K} \ell(i, j) \le \sum_{(i, j) \in K} \frac{\ell_i + \ell_j - 1}{2} + \frac{r}{4}.
    \]
\end{proposition}

\begin{proof}
    Partition $[r] \eqdef \{1, \dots, r\}$ in to index sets $A, B, C, D$ such that 
    \begin{align*}
    A &\eqdef \bigl\{ i \in [r] : \text{$R^i$ is a $+$-path and $\ell_i$ is odd} \bigr\}, \\
    B &\eqdef \bigl\{ i \in [r] : \text{$R^i$ is a $+$-path and $\ell_i$ is even} \bigr\}, \\
    C &\eqdef \bigl\{ i \in [r] : \text{$R^i$ is a $-$-path and $\ell_i$ is odd} \bigr\}, \\
    D &\eqdef \bigl\{ i \in [r] : \text{$R^i$ is a $-$-path and $\ell_i$ is even} \bigr\}. 
    \end{align*}

    For each $(i, j) \in K$, we define $\lambda(i, j) \eqdef \ell_i + \ell_j - 2\ell(i, j)$. We claim that $\lambda(i, j)$ is determined by the parts of $[r]$ that $i$ and $j$ are from, respectively. Here we only go through one of the cases, and the other cases are similar. To be specific, we show that $i, j \in A$ implies $\lambda(i, j) = 2$. In this case, 
    \[
    \ell(i, j) \overset{(*)}{=} \Bigl\lfloor \frac{\ell_i+\ell_j-1}{2}\Bigr\rfloor \overset{(**)}{=} \frac{\ell_i+\ell_j-2}{2} \implies \lambda(i, j) = \ell_i+\ell_j-2\ell(i, j) = 2,  
    \]
    where at $(*)$ we applied that $R^i, R^j$ are of the same sign, and at $(**)$ we used that $\ell_i, \ell_j$ are odd. By working through every other case, we have the table below of values of $\lambda(i, j)$ in all $4 \times 4 = 16$ cases, where the rows and columns stand for the parts that $i$ and $j$ belong to, respectively. 

    \vspace{0.5em}
    \begin{center}
        \begin{tabular}{ |c|c|c|c|c| } 
        \hline
            $\lambda$ & $A$ & $B$ & $C$ & $D$\\ 
        \hline
            $A$ & $2$ & $1$ & $0$ & $1$\\ 
        \hline
            $B$ & $1$ & $2$ & $1$ & $0$\\
        \hline
            $C$ & $0$ & $1$ & $2$ & $1$\\
        \hline
            $D$ & $1$ & $0$ & $1$ & $2$\\ 
        \hline
        \end{tabular}
    \end{center}

    Notice that $|A| + |B| + |C| + |D| = r$.  From the definition of $K$, we deduce that
    \begin{align*}
    |K| &= \binom{|A|}{2} + \binom{|B|}{2} + \binom{|C|}{2} + \binom{|D|}{2} + |A||B| + |B||C| + |C||D| + |D||A| + |A||C| + |B||D| \\
    &= \frac{|A|^2 + |B|^2 + |C|^2 + |D|^2 - r}{2} + |A||B| + |B||C| + |C||D| + |D||A| + |A||C| + |B||D|. 
    \end{align*}
    It then follows that
    \begin{align} \label{eq:lambdabound}
        \sum_{(i, j) \in K} \lambda(i, j) &= \frac{|A|^2 + |B|^2 + |C|^2 + |D|^2 - r}{2} \cdot 2 + \bigl( |A||B| + |B||C| + |C||D| + |D||A| \bigr) \cdot 1 \nonumber \\
        &= |K| + \biggl( \frac{|A|^2 + |B|^2 + |C|^2 + |D|^2 - r}{2} - \bigl( |A||C| + |B||D| \bigr) \biggr) \nonumber \\
        &= |K| + \frac{\bigl(|A|-|C|\bigr)^2 + \bigl(|B|-|D|\bigr)^2 - r}{2} \ge |K| - \frac{r}{2}. 
    \end{align}
    Thus, from the definition of $\lambda$ and \eqref{eq:lambdabound} we obtain
    \[
    \Biggl( \sum_{(i, j) \in K} \frac{\ell_i + \ell_j - 1}{2} + \frac{r}{4} \Biggr) - \Biggl( \sum_{(i, j) \in K} \ell(i, j) \Biggr) = \sum_{(i, j) \in K} \frac{\lambda(i, j)}{2} + \frac{\frac{r}{2} - |K|}{2} \ge 0. \qedhere
    \]
\end{proof}

\medskip
Finally, we are in a position to prove \Cref{thm:fine_covering}. Indeed, 
\begin{align*}
    (2k+1)|S| &= (2k+1) \cdot \sum_{(i, j) \in K} |S_{i, j}| \\
    &\leq (2k+1) \cdot \sum_{(i, j) \in K}
    \ell(i,j) &&\text{by \Cref{prop:Sij_bound}} \\
    &\leq (2k+1) \cdot \Biggl( \sum_{(i, j) \in K} \frac{\ell_i + \ell_j - 1}{2}+\frac{r}{4} \Biggr) &&\text{by \Cref{prop:lij_bound}} \\
    &= (2k+1) \cdot \Biggl( \frac{r-1}{2} \cdot \biggl( \sum_{i=1}^r \ell_i-\frac{r}{2} \biggr) + \frac{r}{4} \Biggr) \\
    &\leq (2k+1) \cdot \Biggl( \frac{r-1}{2} \cdot \biggl(2m(S)-\Bigl(k+\frac12\Bigr)r\biggr) + \frac{r}{4} \Biggr) &&\text{by \Cref{obs:orbit_sum}}. 
\end{align*}
Take the substitution $r' \eqdef (k + \frac{1}{2})r$, and so $r = \frac{2r'}{2k+1}$. This implies that 
\begin{align*}
    (2k+1)|S| &\le (2k+1) \cdot \frac{r-1}{2} \cdot \bigl( 2m(S) - r' \bigr) + \frac{r'}{2} \\
    &= \Bigl( r' - \frac{2k+1}{2} \Bigr) \cdot \bigl( 2m(S) - r' \bigr) + \frac{r'}{2} \\
    &= -(2k+1)m(S) + r' \cdot \bigl( 2m(S) + k + 1 - r' \bigr) \\
    &\overset{(*)}{\le} -(2k+1)m(S) + \Bigl( \frac{2m(S) + k + 1}{2} \Bigr)^2 \\
    &= \Bigl( m(S) - \frac{k}{2} \Bigr)^2 + \frac{2k+1}{4}, 
\end{align*}
where at the step marked with $(*)$ we applied the AM-GM inequality. By rearranging the inequality and taking square roots on both sides, we obtain $m(S) \geq \sqrt{(2k+1)(|S|-\frac14)} + \frac{k}{2}$. 

\section{Generic sets without long \texorpdfstring{$k$}{k}-modal paths} \label{sec:construction}

This section is devoted to the proof of \Cref{thm:modal_setconstruct}. We begin by constructing a family of generic sets $U^{s, t} \subset \R^2$ containing only short $k$-modal paths. Then we can prove \Cref{thm:modal_setconstruct} by carefully choosing $s, t$ as functions of $n$, and delete a small number of points from $U^{s, t}$ to obtain the desired set $S$. It is worth mentioning that our construction is partly inspired by the previous work \cite{bukh_dong_gong_gu_yang_yu}. 

\smallskip

Assume $s,t$ are positive integers. Let $P^1, \dots, P^{s+2t}$ be decreasing paths of lengths 
\[
\underbrace{t, \quad t+1, \quad \dots, \quad 2t-1}_t, \quad \underbrace{2t, \quad \dots, \quad 2t}_{s}, \quad \underbrace{2t-1, \quad \dots, \quad t+1, \quad t}_t
\]
satisfying $p_i \prec_{_{\text{NE}}} p_j$ for any $p_i \in P_i, \, p_j \in P_j$ with $i < j$. Take $U^{s, t} \eqdef P^1 \cup \dots \cup P^{s+2t}$. Then 
\[
|U^{s, t}| = 2\bigl(t + (t+1) + \dots + (2t-1)\bigr) + s \cdot 2t = t(2s+3t-1). 
\]
For instance, \hyperlink{figtwo}{Figure~2}  illustrates a possible configuration of $U^{3, 2}$. 

\centerline{
\begin{tikzpicture}[scale = 0.9]
    \clip (-0.5, -1) rectangle (7.5, 7.5);
    \foreach \x in {0, 1, 2, 3, 4, 5, 6}
        \draw[dash pattern = on 2pt off 2pt] (\x, \x) rectangle (\x+1, \x+1); 
    \foreach \x in {0, 1, 2, 3, 4, 5, 6}
        \filldraw (\x+0.2, \x+0.8) circle (0.05);
    \foreach \x in {0, 1, 2, 3, 4, 5, 6}
        \filldraw (\x+0.8, \x+0.2) circle (0.05);
    \foreach \x in {1, 5}
        \filldraw (\x+0.5, \x+0.5) circle (0.05);
    \foreach \x in {2, 3, 4}
        \filldraw (\x+0.4, \x+0.6) circle (0.05);
    \foreach \x in {2, 3, 4}
        \filldraw (\x+0.6, \x+0.4) circle (0.05);
    \foreach \x in {1, 2, 3, 4, 5, 6, 7}
        \node at (\x-0.775, \x-0.825) {\scriptsize $P^{\x}$};
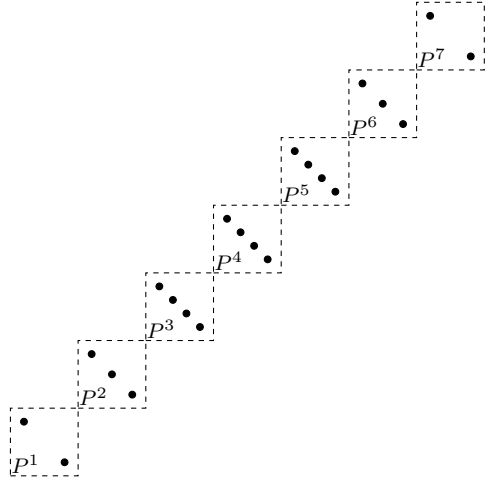
    \node at (3.5, -0.6) {\textbf{\hypertarget{figtwo}{Figure 2:}} One possible drawing of $U^{3, 2}$. };
\end{tikzpicture}
}

The following upper bound on the length of any $k$-modal path subset of $U^{s, t}$ is crucial. 

\begin{proposition} \label{prop:MU_upper}
    Suppose $Q \subseteq U^{s, t}$ is a $k$-modal path. Then $|Q| \le (k+2)t + s - \frac{k}{2}$. 
\end{proposition}

To prove \Cref{prop:MU_upper}, we need some preparations. Suppose $a, b$ are the first and the last point of $Q$, respectively. Find $\alpha$ and $\beta$ such that $a \in P^{\alpha}$ and $b \in P^{\beta}$. Define $Q^{\alpha} \eqdef Q \cap P^{\alpha}, \, Q^{\beta} \eqdef Q \cap P^{\beta}$ and $Q^{\gamma} \eqdef Q \cap P^{\gamma} \, (\gamma = \alpha+1, \dots, \beta-1)$. Call an index $\gamma \in \Gamma^Q \eqdef \{\alpha+1, \dots, \beta-1\}$ as $Q$-\emph{rich} if $|Q \cap P^{\gamma}| \ge 2$, and $Q$-\emph{poor} if $|Q \cap P^{\gamma}| \le 1$. Consider the partition $\Gamma^Q = \Gr \cup \Gp$, where
\[
\Gr \eqdef \bigl\{\gamma \in \Gamma^Q : \text{$\gamma$ is $Q$-rich}\bigr\}, \qquad \Gp \eqdef \bigl\{\gamma \in \Gamma^Q : \text{$\gamma$ is $Q$-poor}\bigr\}. 
\]


\begin{observation} \label{obs:modal_construction}
    If $Q$ is a $(+k)$-path, then $\bigl|\Gr\bigr| \le \lfloor \frac{k}{2} \rfloor$; if $Q$ is a $(-k)$-path, then $\bigl|\Gr\bigr| \le \lfloor \frac{k-1}{2} \rfloor$. 
\end{observation}

\begin{proof}
    The cardinality $\bigl|\Gr\bigr|$ is upper bounded by the number of indices $i \in \{0, 1, \dots, k\}$ such that $Q_i$ is a decreasing section with $|Q_i| \ge 2$. Then the observation follows from counting such sections when $Q$ is a $(+k)$- or $(-k)$-path and $k$ is even or odd, respectively. 
\end{proof}

We say that $\alpha$ is \emph{bad} if the beginning section $Q_0$ of $Q$ is increasing and $|Q^{\alpha}| > 1$. We say that $\beta$ is \emph{bad} if the ending section $Q_k$ of $Q$ is increasing and $|Q^{\beta}| > 1$.

\begin{lemma}
\label{lem:rich_est}
    Denote by $[\mathsf{E}]$ the indicator function of the event $\mathsf{E}$. Then we have 
    \[
    \bigl| \Gr \bigr| \le \lfloor \tfrac{k}{2} \rfloor - [\text{$\alpha$ is bad}] - [\text{$\beta$ is bad}]. 
    \]
\end{lemma}

\begin{proof}
    Observe that $Q_1 \subseteq P^\alpha$ if $\alpha$ is bad and $Q_{k-1} \subseteq P^\beta$ if $\beta$ is bad. We see that there are at most $\lfloor \tfrac{k}{2} \rfloor - [\text{$\alpha$ is bad}] - [\text{$\beta$ is bad}]$ decreasing sections of $Q$ possibly lying in the union $P^{\alpha+1} \cup \dots \cup P^{\beta-1}$. It follows that $\bigl| \Gr \bigr| \le \lfloor \frac{k}{2} \rfloor - [\text{$\alpha$ is bad}] - [\text{$\beta$ is bad}]$. 
\end{proof}

\begin{proof}[Proof of \Cref{prop:MU_upper}]
We prove by a casework on the sign of $Q$ and the parity of $k$. The strategy is to deal with $Q^{\alpha}, \, Q^{\beta}, \, Q^{\gamma} \eqdef Q \cap P^{\gamma} \, (\gamma = \alpha+1, \dots, \beta-1)$ separately. Observe that
\[
|Q^{\alpha}| \le |P^{\alpha}| \le \min\{\alpha+t-1, 2t\}, \qquad |Q^{\beta}| \le |P^{\beta}| \le \min\{s+3t-\beta, 2t\}. 
\]
We shall use one of the estimates on $Q^{\alpha}$ and one of the estimates on $Q^{\beta}$ at each of the $(\spadesuit)$ steps. 
\begin{itemize}
    \item If $Q$ is a $(+k)$-path and $k$ is even, then it follows from the facts $\alpha \ge 1$ and $\beta \le s + 2t$ that $\beta - \alpha - 1 \le s + 2t - 2$. So, from \Cref{lem:rich_est} we deduce that
    \begin{align*}
    |Q| &= |Q^{\alpha}| + |Q^{\beta}| + \sum_{\gamma \in \Gr} |Q^{\gamma}| + \sum_{\gamma \in \Gp} |Q^{\gamma}|\le |Q^{\alpha}| + |Q^{\beta}| + (2t-1) \cdot \bigl| \Gr \bigr| + \bigl| \Gamma^Q \bigr| \\
    &\leq |Q^{\alpha}| + |Q^{\beta}| + (2t-1) \cdot \bigl(\tfrac{k}{2} - \bigl[ |Q^\alpha| > 1 \bigr] - \bigl[ |Q^\beta| > 1 \bigr]\bigr) + (\beta - \alpha - 1) \\
    &\overset{(\spadesuit)}{\leq} 1 + 1 + (2t-1) \cdot \tfrac{k}{2} + (s+2t-2) = (k+2)t + s - \tfrac{k}{2}. 
    \end{align*}
    
    \item If $Q$ is a $(-k)$-path and $k$ is even, then from \Cref{obs:modal_construction} we deduce that
    \begin{align*}
    |Q| &= |Q^{\alpha}| + |Q^{\beta}| + \sum_{\gamma \in \Gr} |Q^{\gamma}| + \sum_{\gamma \in \Gp} |Q^{\gamma}| \le |Q^{\alpha}| + |Q^{\beta}| + (2t-1) \cdot \bigl|\Gr\bigr| + \bigl| \Gamma^Q \bigr| \\
    &\overset{(\spadesuit)}\le (\alpha+t-1) + (s+3t-\beta) + (2t-1) \cdot \tfrac{k-2}{2} + (\beta-\alpha-1) = (k+2)t + s - \tfrac{k}{2} - 1. 
    \end{align*}
    
    \item If $Q$ is a $(+k)$-path and $k$ is odd, then from \Cref{lem:rich_est} and the fact $\alpha \ge 1$ we deduce that
    \begin{align*}
    |Q| &= |Q^{\alpha}| + |Q^{\beta}| + \sum_{\gamma \in \Gr} |Q^{\gamma}| + \sum_{\gamma \in \Gp} |Q^{\gamma}| \le |Q^{\alpha}| + |Q^{\beta}| + (2t-1) \cdot \bigl|\Gr\bigr| + \bigl| \Gamma^Q \bigr| \\
    &\leq |Q^{\alpha}| + |Q^{\beta}| + (2t-1) \cdot \bigl( \tfrac{k-1}{2} - \bigl[|Q^\alpha| > 1\bigr] \bigr) + (\beta - \alpha - 1) \\
    &\overset{(\spadesuit)}\le 1 + (s+3t-\beta) + \tfrac{k-1}{2} \cdot (2t-1) + (\beta-1-1) = (k+2)t + s - \tfrac{k}{2} - \tfrac{1}{2}. 
    \end{align*}
    
    \item If $Q$ is a $(-k)$-path and $k$ is odd, then from \Cref{lem:rich_est} and the fact $\beta \le s + 2t$ we deduce that
    \begin{align*}
    |Q| &= |Q^{\alpha}| + |Q^{\beta}| + \sum_{\gamma \in \Gr} |Q^{\gamma}| + \sum_{\gamma \in \Gp} |Q^{\gamma}| \le |Q^{\alpha}| + |Q^{\beta}| + (2t-1) \cdot \bigl|\Gr\bigr| + \bigl| \Gamma^Q \bigr| \\
    &\leq |Q^{\alpha}| + |Q^{\beta}| + (2t-1) \cdot \bigl(\tfrac{k-1}{2} - \bigl[ |Q^\beta| > 1 \bigr]\bigr) + (\beta - \alpha - 1) \\
    &\overset{(\spadesuit)} \le (\alpha+t-1) + 1 + \tfrac{k-1}{2} \cdot (2t-1) + (s+2t-\alpha-1) = (k+2)t + s - \tfrac{k}{2} - \tfrac{1}{2}. 
    \end{align*}
\end{itemize}
By combining all four cases above, we conclude that $|Q| \le (k+2)t + s - \frac{k}{2}$. 
\end{proof}

\begin{proof} [Proof of \Cref{thm:modal_setconstruct}]
    For given integers $k \ge 1$ and $n \ge 10k^3$, we choose parameters
    \[
    x \eqdef \Bigl\lceil \sqrt{(2k+1)(n-\tfrac{1}{4})} - \tfrac{k}{2} \Bigr\rceil, \qquad y \eqdef x + 1 + \lceil \tfrac{k}{2} \rceil, \qquad t \eqdef \langle \tfrac{y}{2k+1} \rangle, \qquad s \eqdef y - (k+2)t, 
    \]
    where $\langle \alpha \rangle \eqdef \lfloor \alpha + \frac{1}{2} \rfloor$ denotes the closest integer to $\alpha$. For any generic $T \subset \R^2$, denote by $M(T)$ the maximum length of a $k$-modal path subset of $T$. Write $U \eqdef U^{s, t}$ for brevity. 

    We first show that $M(U) \le x + 1$. Indeed, it follows from \Cref{prop:MU_upper} that
    \[
    M(U) \le (k+2)t + s - \bigl\lceil \tfrac{k}{2} \bigr\rceil = y - \bigl\lceil \tfrac{k}{2} \bigr\rceil = x+1. 
    \]
    
    We then show that $|U| \ge n$. From the definition we deduce that $\frac{y-k}{2k+1} \le t \le \frac{y+k}{2k+1}$. So, 
    \begin{align*}
    |U| = t(2s+3t-1) &= t\bigl(2y-(2k+1)t+1\bigr) = -(2k+1) \bigl( t - \tfrac{y-\frac12}{2k+1} \bigr)^2 + \tfrac{(y-\frac12)^2}{2k+1} \\
    &\ge -(2k+1) \bigl( \tfrac{k+\frac12}{2k+1} \bigr)^2 + \tfrac{(y-\frac12)^2}{2k+1} = \tfrac{(y-\frac12)^2}{2k+1} - \tfrac{2k+1}{4}. 
    \end{align*}
    It follows that
    \[
    |U| \ge n \impliedby \tfrac{(y-\frac12)^2}{2k+1} - \tfrac{2k+1}{4} \ge n \iff (2y-1)^2 \ge 4(2k+1)n + (2k+1)^2. 
    \]
    By $x \ge \sqrt{(2k+1)(n-\frac14)} - \frac{k}{2}$ and $y \ge x + 1 + \frac{k}{2}$ we have $2y-1 \ge \sqrt{(2k+1)(4n-1)} + 1$, and so
    \[
    (2y-1)^2 = \bigl( \sqrt{(2k+1)(4n-1)} + 1 \bigr)^2 \ge 4(2k+1)n + (2k+1)^2, 
    \]
    thanks to the assumption $n \ge 10k^3$. This implies that $|U| \ge n$. 

    Finally, since $|U| \ge n$, we fix choose an arbitrary $n$-point subset of $U$ as $S$. It follows that 
    \[
    M(S) \le M(U) \le x+1 = \Bigl\lceil \sqrt{(2k+1)(n-\tfrac{1}{4})} - \tfrac{k}{2} \Bigr\rceil + 1. 
    \]
    So, the set $S$ constructed above satisfies \Cref{thm:modal_setconstruct}, and the proof is complete. 
\end{proof}

\section*{Acknowledgements}

I am grateful to Zichao Dong for suggesting this problem, for talking with me the previous work \cite{bukh_dong_gong_gu_yang_yu} (the undergraduate research was supervised by Boris Bukh and Zichao Dong), and for many helpful discussions. I would also like to thank Professor Chunwei Song for providing valuable suggestions on improving the writing of this paper. 

\bibliographystyle{plain}
\bibliography{M-seq-final}

\appendix

\section{Proof of \texorpdfstring{\Cref{lem:inj}}{Lemma~\ref{lem:inj}}} \label{sec:proof_inj}

For any $P \in \cC_+ \cup \cC_-$, let $\gamma(P) \in \Z_{\ge 0}$ be the maximum index $j$ such that $\sN_P$ sends the first through the $j$-th point of $P$ to a same $Q$, which is $\varphi(P)$ by definition. Here we set $\gamma(P)$ to be $0$ if $i(P) = k$. Write $\vH_i \eqdef \sum\limits_{\substack{P \in \cC_+ \cup \cC_- \\ d(P) = i}} \gamma(P)^2$ for $i = 0, 1, \dots$ and associate with $\cC_+, \cC_-$ the infinite tuple
\[
\vH \eqdef (\vH_0, \vH_1, \dots) \in (\Z_{\ge 0})^{(\Z_{\ge 0})}. 
\]
We claim that $\vH$ takes finitely many possible values. Indeed, from our assumption at the beginning of \Cref{sec:prop} we obtain $|\cC_+ \cup \cC_-| \le 2m(S)$. Also, we have $\gamma(P) \le |S|$ and $d(P) \le |\cC_+ \cup \cC_-| \le 2m(S)$ (because the action of $\varphi$ is acyclic). It follows that $\vH_i$ is upper bounded by $2|S|^2 \cdot m(S)$. 

Among all possible pairs of $\cC_+, \cC_-$ with $\max\bigl\{|\cC_+|, |\cC_-|\bigr\} \le m(S)$, we choose one so that $\vH$ is maximized under the canonical lexicographical order on $(\Z_{\ge 0})^{(\Z_{\ge 0})}$. Our aim is to establish that $\varphi$, associated with this choice of $\cC_+, \cC_-$, is injective. Notice that $\sN_P, \cD, \varphi, d, \gamma, \vH$ depend on $\cC_+, \cC_-$. 

Assume to the contrary that $P = \varphi(Q) = \varphi(R)$ for some $Q \neq R$. This implies that $Q, R$ are of the same sign, for $P$ and $\varphi(P)$ are always from different ones of $\cC_+, \cC_-$. That is, the monotonicity of $Q_i$ and $R_i$ are the same for each $i$. Keep in mind that even if $|Q_i| \le 1$ for some $i$, we still endow a fixed monotonicity to $Q_i$ as $Q$ is from a certain one of $\cC_+$ and $\cC_-$. 

Since $\varphi(Q) = \varphi(R)$ implies $d(Q) = d(R) = d$, we may choose $Q$ and $R$ so that $d$ is minimized. Let $q$ be the $\gamma(Q)$-th point of $Q$ and $r$ be the $\gamma(R)$-th point of $R$. Assume without loss of generality that $x(q) \ge x(r)$. We claim that $x(q) > x(r)$. If not, then $x(q) = x(r)$, and so $q = r = p$. Suppose $p \in P_{j+1}$. Then from $\varphi(Q) = \varphi(R) = P$ we deduce that $q \in Q_j, \, r \in R_j$. However, it follows from $\cC_+, \cC_-$ are fine coverings and $Q, R$ are of the same sign that $Q_j, R_j$ are disjoint, a contradiction. 

Informally speaking, we are going to throw a couple of points of $R$ into $Q$, and hence produce another pair of fine coverings $\cC_+', \cC_-'$ with $\vH' \succ \vH$, which is a contradiction. 

Suppose $A \subseteq R$ is the set of the $1$-st through the $\gamma(R)$-th point of $R$, and $B \subseteq Q$ is the set of the $1$-st through the $\gamma(Q)$-th point of $Q$. Observe from the definition of $\gamma$ that $R_i \cap A \subseteq P_{i+1}$ and $Q_i \cap B \subseteq P_{i+1}$ hold for all $i = 0, 1, \dots, k$. We obtain $\cC_+'\cup\cC_-'$ from $\cC_+ \cup \cC_-$ by replacing $Q, R$ with $Q', R'$ and keep every other $k$-modal path intact, where for $i = 0, 1, \dots, k$, 
\begin{itemize}
    \item the $i$-th sections of $Q', R'$ are defined as $Q_i' \eqdef Q_i \cup (R_i\cap A)$ and $R_i' \eqdef R_i \setminus A$, respectively. 
\end{itemize}
Then $Q_i' \cup R_i' = Q_i \cup R_i$ as multisets. Obviously, $R'$ is a $k$-modal path of the same sign as $R$. 

We claim that $Q'$ is a $k$-modal path of the same sign as $Q$. To see this, we need to show that 
\begin{enumerate}[label=(Q\arabic*), ref=(Q\arabic*)]
	\item \label{Qiprime_single} $Q_i'$ shares the same monotonicity with $Q_i$ for each $i = 0, 1, \dots, k$, and
	\item \label{Qiprime_pair} $\max\limits_{p \in Q_s'} x(p) < \min\limits_{p \in Q_t'} x(p)$ for any pair of indices $s, t$ with $0 \le s < t \le k$. 
\end{enumerate}
For \ref{Qiprime_single}, since $Q_i, R_i$ are of the same monotonicity with $P_{i+1}$, it suffices to show that $Q_i'$ is monotone. Choose index $t$ such that $q \in Q_t$. We prove \ref{Qiprime_single} through the following casework on $i$: 
\begin{itemize}
    \item If $i < t$, then $Q_i' = Q_i \cup (R_i \cap A) = (Q_i \cap B) \cup (R_i \cap A) \subseteq P_{i+1}$. To see the inclusion step, the definition of $B$ implies that $\sN_Q(b) = P$ for every $b \in B$, and so $Q_i \cap B \subseteq P_{i+1}$. Similarly, we have $R_i \cap A \subseteq P_{i+1}$. Since $P$ is a $k$-modal path, $Q_i' \subseteq P_{i+1}$ implies that $Q_i'$ is monotone. 
    \item If $i > t$, then from $x(q) > x(r)$ we deduce that $Q_i' = Q_i$, and so $Q_i'$ is monotone. 
    \item If $i = t$, then we consider $X \eqdef \bigl\{a \in Q_i' : x(a) \leq x(q)\bigr\}, \, Y \eqdef \bigl\{a \in Q_i' : x(a) \geq x(q)\bigr\}$ separately. 
    \begin{itemize}
        \item By applying $x(q) > x(r)$ at $(*)$, we may expand $X$ into a union of two sets as follows: 
        \begin{align*}
        X &= \bigl\{a \in Q_i \cup (R_i \cap A) : x(a) \le x(q)\bigr\} \\
        &= \bigl\{a \in Q_i : x(a) \le x(q)\bigr\} \cup \bigl\{a \in R_i \cap A : x(a) \le x(q)\bigr\} \\
        &\overset{(*)}{=} \bigl\{a \in Q_i : x(a) \le x(q)\bigr\} \cup \bigl\{a \in R_i : x(a) \le x(r)\bigr\}. 
        \end{align*}
        The definition of $\gamma$ tells that both sets above are contained in $P_{i+1}$, and so $X \subseteq P_{i+1}$. 
        \item Again, from the fact $x(q) > x(r)$ we deduce that $Y \subseteq Q_i$. 
    \end{itemize}
    So, $Q_i' = X \cup Y$ is monotone, thanks to $Q_i, P_{i+1}$ are of the same monotonicity and $q \in X \cap Y$. 
\end{itemize}
Recall the conventions $\max \varnothing = -\infty$ and $\min \varnothing = +\infty$. We will prove \ref{Qiprime_pair} by showing $x(b) < x(c)$ for any $(b, c) \in Q_i' \times Q_j'$ with $0 \le i < j \le k$ through the following casework on $i, j$: 
\begin{itemize}
    \item If $0 \leq i < j \leq t$, then we have shown $Q_i' \subseteq P_{i+1}$ in the proof of \ref{Qiprime_single}. Consider the first point $a$ of $Q_j'$ (there is nothing to prove if $Q_j' = \varnothing$). Then $x(a) \leq x(q)$. Recall that $Q_j' = Q_j \cup (R_j \cap A)$. 
    \begin{itemize}
        \item If $a \in Q_j$, then $a \in B$ thanks to $x(a) \le x(q)$, and so $a \in P_{j+1}$. 
        \item If $a \in R_j \cap A$, then from $\varphi(R) = P$ we deduce that $a \in P_{j+1}$. 
    \end{itemize}
    Since $P$ is a $k$-modal path and $i+1 < j+1$, we obtain $x(b) < x(a) \le x(c)$. 
    \item For $t \leq i < j \leq k$, then we have shown $Q_j' = Q_j$ in the proof of \ref{Qiprime_single}. Consider the last point $a$ of $Q_i'$ (there is nothing to prove if $Q_i' = \varnothing$). Then $x(a) \geq x(q) > x(r)$, and so $a \notin A$, hence $a \notin R_j \cap A$. Thus, $a \in Q_i$. Since $Q$ is a $k$-modal path and $i < j$, we obtain $x(b) \le x(a) < x(c)$. 
    \item If $0 \le i < t < j \le k$, then by considering $(\widetilde{i}, \widetilde{j}) = (i, t)$ with $0 < \widetilde{i} < \widetilde{j} \le t$ and $(\widetilde{i}, \widetilde{j}) = (t, j)$ with $t \le \widetilde{i} < \widetilde{j} \le k$ in what we just proved, it follows that $x(b) < x(q) < x(c)$. 
\end{itemize}

By combining the properties of $Q'$ and $R'$ proved above, we conclude that $\cC_+'$ and $\cC_-'$ are a $(+k)$- and a $(-k)$-covering of $S$, respectively. Define $\sN_{P'}', \cD', \varphi', d', \gamma', \vH'$ with respect to $\cC_+'$ and $\cC_-'$ in the same way as $\sN_P, \cD, \varphi, d, \gamma, \vH$ with respect to $\cC_+$ and $\cC_-$. For any $L \in \cC_+ \cup \cC_-$, write $L' \eqdef Q', R'$ when $L = Q, R$ (respectively) and $L' \eqdef L$ otherwise. 

\begin{proposition} \label{prop:phi_commute}
    We have $\varphi'(L') = \varphi(L)'$ for all $L \in \cD$ with $L \neq R$ and $\varphi(L) \neq R$. 
\end{proposition}

\begin{proof}
    An important observation is that $\varphi(L) = \sN_L\bigl(s(L)\bigr)$ is uniquely determined by $\bigl(s(L), i(L)\bigr)$ in $\cC_+ \cup \cC_-$. Similarly, $\varphi'(L') = \sN'_{L'}\bigl(s(L')\bigr)$ is uniquely determined by $\bigl( s(L'), i(L') \bigr)$ in $\cC_+' \cup \cC_-'$. 
    \begin{itemize}
        \item If $L = Q$, then $\bigl(s(L'), i(L')\bigr) = \bigl(s(Q'), i(Q')\bigr) \in \bigl\{ \bigl(s(Q),i(Q)\bigr), \bigl(s(R), i(R)\bigr) \bigr\}$ since $Q'$ consists of $Q$ and some part of $R$. It follows from $P = \varphi(R) = \varphi(Q) \notin \{Q, R\}$ that 
        \[
        \varphi'(Q') = \sN_{Q'}'\bigl(s(Q')\bigr) \in \bigl\{\sN_Q\bigl(s(Q)\bigr), \sN_R\bigl(s(R)\bigr)\bigr\} = \bigl\{\varphi(Q), \varphi(R)\bigr\} = \{P\} = \{P'\} = \{\varphi(Q)'\}. 
        \]
        \item If $L \neq Q$, then from $L \neq R$ we deduce that $L' = L$, and hence $\bigl(s(L'), i(L')\bigr) = \bigl(s(L), i(L)\bigr)$. Let $\varphi(L) \eqdef \widetilde{P}$. Then $(\widetilde{P})_i \subseteq (\widetilde{P}')_i$, where the strict inclusion could happen only when $\widetilde{P} = Q$. It follows that $\widetilde{P}'$ contains $s(L')$, and hence $\varphi'(L') = \widetilde{P'} = \varphi(L)'$. 
    \end{itemize}
    By combining the cases above, we conclude that the proposition holds. 
\end{proof}


Recall that $d(Q) = d(R) = d$. For $i = 0, 1, \dots, d$, write 
\[
\cL_i \eqdef \{L \in \cC_+ \cup \cC_- : d(L) = i\}, \qquad \cL'_i \eqdef \{L '\in \cC'_+ \cup \cC'_- : d'(L') = i\}. 
\]
Denote $\cL_{\le d} \eqdef \bigcup\limits_{i=0}^d \cL_i$ and $\cL'_{\le d} \eqdef \bigcup\limits_{i=0}^d \cL'_i$. Then we have for every $L \in \cL_{\le d} \setminus \{Q, R\}$ that $L' = L$. We claim first for each $L \in \cL_{\le d}$ that $\bigl(d'(L'), \gamma'(L')\bigr) = \bigl(d(L), \gamma(L)\bigr)$, and second that $d'(Q') = d$. 
\begin{itemize}
    \item Firstly, $\varphi^{(i)}(L) \neq R$ for $L \in \cL_{\le d} \setminus \{Q, R\}$ because otherwise the depth of $L$ would be bigger than $d$. So, by \Cref{prop:phi_commute}, we see that $\varphi^{(i)}(L)' = \varphi'^{(i)}(L')$, and hence $d'(L') = d(L)$. From $L = L'$ and $\varphi(L) = \varphi(L)'$ we obtain $\gamma'(L') = \gamma(L)$. 
    \item Secondly, from $\varphi(Q) = P$ we deduce that $d(P) = d(Q) - 1$. Since \Cref{prop:phi_commute} implies that $\varphi'(Q') = P'$, we obtain from the first claim that $d'(Q') = d'(P') + 1 = d(P) + 1 = d$. 
\end{itemize}

We next show that $\gamma'(Q')^2 > \gamma(Q)^2 + \gamma(R)^2$. It suffices to prove $\gamma'(Q') \ge \gamma(Q) + \gamma(R)$. (In fact, the equality holds, but we shall not need this stronger fact.) Indeed, we are to show that every point of $Q_i'$ which is to the left of $q$ contributes to $\gamma'(Q')$. That is, $\sN'_{Q'}\bigl(s(Q')\bigr) = \dots = \sN'_{Q'}(q)$. Recall that $Q, R$ are part of a fine covering and so $Q_i \cap R_i = \varnothing$. This will then imply that $\gamma'(Q') \ge \gamma(Q) + \gamma(R)$ as $x(r) < x(q)$ and every point of $R$ that contributes to $R$ are moved into $Q'$. Fix an index $i$ and consider any point $a \in (Q_i \cap B) \cup (R_i \cap A) \subseteq Q_i'$ (hence $a$ is to the left of $q$). Then $a\in P_{i+1} = P_{i+1}'$ and $\sN_{Q'}'(a) = P$. This implies that $\sN'_{Q'}\bigl(s(Q')\bigr) = \dots = \sN'_{Q'}(q) = P$, and so $\gamma'(Q')\geq \gamma(Q) + \gamma(R)$. 

\smallskip

Finally, we are ready to prove $\vH_i \le \vH_i'$ for $i = 0, 1, \dots, d-1$ and $\vH_d < \vH_d'$. Indeed, for any $L \in \cL_{\le d} \setminus \{Q, R\}$, we have $d'(L') = d(L)$ and $\gamma'(L') = \gamma(L)$. Therefore, by $\gamma'(Q')^2 > \gamma(Q)^2 + \gamma(R)^2$, 
\begin{align*}
    \vH_i = \sum_{L \in \cL_i} \gamma(L)^2 &= \sum_{L' \in \cL'_d} \gamma'(L')^2 = \vH_i' \quad (i = 0, 1, \dots, d-1), \\
    \vH_d = \sum_{L \in \cL_d} \gamma(L)^2 &= \biggl( \sum_{L \in \cL_d \setminus \{Q, R\}} \gamma(L)^2 \biggr) + \bigl(\gamma(Q)^2 + \gamma(R)^2\bigr) \\
    &< \biggl(\sum_{L' \in \cL'_d \setminus \{Q', R'\}} \gamma'(L')^2 \biggr) + \gamma'(Q')^2 \\
    &\leq \sum_{L' \in \cL'_d} \gamma'(L')^2 = \vH_d'.
\end{align*}
Thus, we obtain the desired contradiction $\vH' \succ \vH$. The proof of \Cref{lem:inj} is complete. 

\end{document}